\documentstyle[11pt,amstex,amscd,graphicx]{amsart}

\begin{document}

\newtheorem{theorem}{Theorem}[section]
\newtheorem{prop}[theorem]{Proposition}
\newtheorem{lemma}[theorem]{Lemma}
\newtheorem{cor}[theorem]{Corollary}
\newtheorem{defn}[theorem]{Definition}
\newtheorem{conj}[theorem]{Conjecture}
\newtheorem{claim}[theorem]{Claim}
\newtheorem{rmk}[theorem]{Remark}

\newcommand{\bdy}{\partial}
\newcommand{\C}{{\mathbb C}}
\newcommand{\integers}{{\mathbb Z}}
\newcommand{\natls}{{\mathbb N}}
\newcommand{\ratls}{{\mathbb Q}}
\newcommand{\reals}{{\mathbb R}}
\newcommand\Proj{{\mathbf P}}
\newcommand\U{{\mathbf U}}
\newcommand\Hyp{{\mathbf H}}
\newcommand\D{{\mathbf D}}
\newcommand{\T}{{\mathbf T}}
\newcommand{\I}{{\mathbf I}}
\newcommand\Z{{\mathbb Z}}
\newcommand\R{{\mathbb R}}
\newcommand\Q{{\mathbb Q}}
\newcommand\E{{\mathbb E}}
\newcommand{\proj}{{\mathbb P}}
\newcommand{\lhp}{{\mathbb L}}
\newcommand\AAA{{\mathcal A}}
\newcommand\BB{{\mathcal B}}
\newcommand\CC{{\mathcal C}}
\newcommand\DD{{\mathcal D}}
\newcommand\EE{{\mathcal E}}
\newcommand\FF{{\mathcal F}}
\newcommand\GG{{\mathcal G}}
\newcommand\HH{{\mathcal H}}
\newcommand\II{{\mathcal I}}
\newcommand\JJ{{\mathcal J}}
\newcommand\KK{{\mathcal K}}
\newcommand\LL{{\mathcal L}}
\newcommand\MM{{\mathcal M}}
\newcommand\NN{{\mathcal N}}
\newcommand\OO{{\mathcal O}}
\newcommand\PP{{\mathcal P}}
\newcommand\QQ{{\mathcal Q}}
\newcommand\RR{{\mathcal R}}
\newcommand\SSS{{\mathcal S}}
\newcommand\TT{{\mathcal T}}
\newcommand\UU{{\mathcal U}}
\newcommand\VV{{\mathcal V}}
\newcommand\WW{{\mathcal W}}
\newcommand\XX{{\mathcal X}}
\newcommand\YY{{\mathcal Y}}
\newcommand\ZZ{{\mathcal Z}}
\newcommand\CH{{\CC\HH}}
\newcommand\MF{{\MM\FF}}
\newcommand\PMF{{\PP\kern-2pt\MM\FF}}
\newcommand\ML{{\MM\LL}}
\newcommand\PML{{\PP\kern-2pt\MM\LL}}
\newcommand\GL{{\GG\LL}}
\newcommand\Pol{{\mathcal P}}
\newcommand\half{{\textstyle{\frac12}}}
\newcommand\Half{{\frac12}}
\newcommand\ep{\epsilon}
\newcommand\hhat{\widehat}
\newcommand\til{\widetilde}
\newcommand\gesim{\succ}
\newcommand\lesim{\prec}
\newcommand\simle{\lesim}
\newcommand\simge{\gesim}
\newcommand{\simmult}{\asymp}
\newcommand{\simadd}{\mathrel{\overset{\text{\tiny $+$}}{\sim}}}
\newcommand{\sm}{\setminus}
\newcommand{\pair}[1]{\langle #1\rangle}
\newcommand{\tprec}{\prec_t}
\newcommand{\fprec}{\prec_f}
\newcommand{\bprec}{\prec_b}
\newcommand{\pprec}{\prec_p}
\newcommand{\ppreceq}{\preceq_p}
\newcommand{\sprec}{\prec_s}
\newcommand{\cpreceq}{\preceq_c}
\newcommand{\cprec}{\prec_c}
\newcommand{\topprec}{\prec_{\rm top}}
\newcommand{\Topprec}{\prec_{\rm TOP}}
\newcommand{\fsub}{\mathrel{\scriptstyle\searrow}}
\newcommand{\bsub}{\mathrel{\scriptstyle\swarrow}}
\newcommand{\fsubd}{\mathrel{{\scriptstyle\searrow}\kern-1ex^d\kern0.5ex}}
\newcommand{\bsubd}{\mathrel{{\scriptstyle\swarrow}\kern-1.6ex^d\kern0.8ex}}
\newcommand{\fsubeq}{\mathrel{\raise-.7ex\hbox{$\overset{\searrow}{=}$}}}
\newcommand{\bsubeq}{\mathrel{\raise-.7ex\hbox{$\overset{\swarrow}{=}$}}}
\newcommand{\bbar}{\overline}
\newcommand{\UML}{\operatorname{\UU\MM\LL}}
\newcommand{\EL}{\mathcal{EL}}
\newcommand\A{\mathbf A}

\title{Relative Rigidity, Quasiconvexity and C-Complexes}

\author{Mahan Mj}
\address{RKM Vivekananda University, Belur Math, WB-711 202, India}

\date{}

\begin{abstract} 
We introduce and study the notion of relative rigidity for pairs $(X,\JJ)$ where \\
1) $X$ is a hyperbolic metric space and $\JJ$ a collection of quasiconvex sets \\
2) $X$ is a relatively hyperbolic group and $\JJ$ the collection of parabolics \\
3) $X$ is a higher rank symmetric space and $\JJ$ an equivariant collection of maximal flats \\ 
Relative rigidity can roughly be described as upgrading  a uniformly proper map between two such $\JJ$'s to a quasi-isometry between the corresponding $X$'s.  A related notion is that of a $C$-complex which is the adaptation of a Tits complex to this context. We prove the relative rigidity of the collection of pairs $(X, \JJ)$ as above. This generalises a result of Schwarz for symmetric patterns of geodesics in hyperbolic space. We show that a uniformly proper map induces an isomorphism of the corresponding $C$-complexes. We also give a couple of characterizations of quasiconvexity
of subgroups of hyperbolic groups on the way.

\begin{center}
AMS subject classification =   20F67(Primary), 22E40,   57M50(Secondary)
\end{center}

\end{abstract}

\maketitle

\tableofcontents

\section{Introduction}

\subsection{Relative Rigidity and Statement of Results}
In this paper, we study a rigidity phenomenon within the framework of coarse geometry. We call it {\em relative rigidity}.  Much of the work on quasi-isometric rigidity (e.g. Farb-Schwartz \cite{farb-schwarz} Kleiner-Leeb \cite{kleiner-leeb} Eskin-Farb \cite{eskin-farb-jams} and Mosher-Sageev-Whyte \cite{msw1} \cite{msw2} ) contains a crucial step showing that a self quasi-isometry of a space $X$ coarsely preserves a family $\JJ$ of distinguished subsets of $X$. The family $\JJ$ again has a coarse intersection pattern that may be combinatorially coded and these proofs of quasi-isometric rigidity often show that the intersection pattern is preserved by a quasi-isometry. In this note, we investigate a sort of a converse to this: \\
{\it When does a uniformly proper map between two families $\JJ_1$ and $\JJ_2$  come from a quasi-isometry $\phi$ between $X_1$ and $X_2$? Does such a map preserve intersection patterns?} \\
We show that the answer is affirmative when \\
\begin{enumerate}
\item $X_i$'s are (Cayley graphs of) hyperbolic groups and $\JJ_i$'s correspond to cosets of a quasiconvex subgroup\\
\item $X_i$'s are (Cayley graphs of) relatively hyperbolic groups and $\JJ_i$'s correspond to parabolic subgroups\\
\item $X_i$'s are symmetric spaces of non-positive curvature and $\JJ_i$'s correspond to lifts of a maximal torus in a compact locally symmetric space modeled on $X_i$. \\
\end{enumerate}

 If in addition one can show that a quasi-isometry preserving intersection patterns is close to an isometry, we would be able to conclude that a uniformly proper map between the $\JJ_i$'s is induced by an isometry. This latter phenomenon has been investigated by Mosher, Sageev and Whyte \cite{msw2} and has been termed {\it pattern rigidity}. 
Thus, in a sense,
 the notion of {\it relative rigidity}
 complements that of {\it pattern rigidity}. 

Some further examples where a family of distinguished subsets of a space 
 and the resulting (combinatorial) configuration
 yields information about the ambient space are: \\
1) Collection of flats in a symmetric space of higher rank (Mostow \cite{mostow-book})\\
2) Collection of maximal abelian subgroups of the mapping class group (Behrstock-Drutu-Mosher \cite{behrstock-drutu-mosher} ) \\
3) Collection of hyperbolic spaces in the Cayley complex of the Baumslag-Solitar groups (Farb-Mosher \cite{farb-mosher-bs1} , \cite{farb-mosher-bs2} ; 
see also \cite{farb-mosher-acta} ) \\
4) Quasi-isometric rigidity of sufficiently complicated patterns of flats in the universal cover of a Haken 3 manifold (Kapovich-Leeb \cite{kap-leeb} )\\
5) We were most influenced by
a beautiful result of Schwarz \cite{schwarz-inv} which shows that a uniformly proper map from a symmetric pattern of geodesics in $\Hyp^n$ to another symmetric pattern of geodesics in $\Hyp^n$ (for $n > 2$) is induced by an isometry. Again as in Mostow, there are two parts to this. A first step  is to construct a quasi-isometry of $\Hyp^n$ inducing the given pairing. Schwarz terms this {\em ambient extension}. The second is to construct an isometry.

\smallskip

Let us look at a general form of the situation that Schwarz considers. $(X_1, d_1), (X_2,d_2)$ are metric spaces. Let $\JJ_1, \JJ_2$ be collections of closed subsets of $X_1, X_2$ respectively. Then $d_i$ induces a {\em pseudo-metric} (which, by abuse of notation, we continue to refer to as $d_i$) on $\JJ_i$ for $i = 1, 2$. This is just the ordinary (not Hausdorff) distance between closed subsets of a metric space. In \cite{schwarz-inv}, $X_1 = X_2 = \Hyp^n$, and $\JJ_i$ are lifts (to the universal cover) of  finite collections of closed geodesics in two hyperbolic manifolds. 

Also, the hypothesis in Schwarz's paper \cite{schwarz-inv} is the existence of a {\em uniformly proper} map $\phi$ between symmetric patterns of geodesics $\JJ_1$ and $\JJ_2$. 
A uniformly proper map may be thought of as an isomorphism in the so-called {\em coarse category } in the sense of John Roe \cite{roe-cbms}. Thus, we can re-interpret the first step of Schwarz's result as saying that an isomorphism $\phi$ in the  coarse category between $\JJ_i$ implies the existence of a quasi-isometry from $\Hyp^n$ to itself inducing $\phi$. In the language of \cite{schwarz-inv}, {\em uniformly proper pairings come from ambient extensions.} 

In Mostow's proof of rigidity for higher rank symmetric spaces, he obtains in a crucial step,
 an isomorphism of Tits complexes \cite{mostow-book}.  We would like to associate to a pair $(X, \JJ)$  some such complex just as a Tits complex is associated to a higher rank locally symmetric space and its collection of maximal parablic subgroups. We propose the notion of a $C$-complex in this paper as the appropriate generalization of a Tits complex to coarse geometry. Then what we would hope for (as a conclusion) is an isomorphism of these $C$-complexes. This transition from the existence of a uniformly proper map between $\JJ_i$'s to the existence of a 
a quasi-isometry between $X_i$'s inducing an isomorphism of $C$-complexes is what we term {\bf relative rigidity}. Schwarz proves the relative rigidity of pairs $(X,\JJ)$ where $X$ is hyperbolic space and $\JJ$ a symmetric collection of geodesics. Much of what he does extends to the case where $X$ is a higher rank symmetric space and $\JJ$ a symmetric collection of maximal periodic flats or a symmmetric collection of maximal parabolic subgroups in a non-uniform lattice.

\smallskip

The main point of this paper is illustrated first in the context of
 relative rigidity of  the category of pairs $(\Gamma , \JJ )$, where $\Gamma$ is (the Cayley graph of) a hyperbolic group, and $\JJ$ the set of cosets of a quasiconvex subgroup. {\bf Throughout this paper we shall assume that the quasiconvex subgroups are of infinite index in the big groups.}

Note that the upgrading of a uniformly proper map between $\JJ$'s to a quasi-isometry between the $\Gamma$'s is the most we can hope for in  light of the fact that the Cayley graph of a finitely generated group is only determined up to quasi-isometry. (See  Paulin \cite{paulin-bdy} for a proof of this fact.)

We start with a pair of hyperbolic groups $G_1, G_2$ with Cayley graphs $\Gamma_1, \Gamma_2$, and quasiconvex subgroups $H_1, H_2$. Let $\Lambda_1$, $\Lambda_2$ be the limit sets of $H_1, H_2$ in $\partial G_1, \partial G_2$ respectively. For convenience we consider the collection $\JJ_i$ of translates of $J_i$ the join of $\Lambda_i$ in $\Gamma_i$ rather than cosets of $H_i$. Recall that the join of $\Lambda_i$ is the union of bi-infinite geodesics in $\Gamma_i$ with end-points in $\Lambda_i$. This is a uniformly quasiconvex set and lies at a bounded Hausdorff distance from the Cayley graph of the subgroup $H_i$ (Since $H$ has finite index in its commensurator, only finitely many cosets of $H$ are at
a finite Hausdorff distance from it. Since $J_i$ is at a bounded Hausdorff distance from $H_i$ the same is true for
elements of $\JJ_i$.) The main theorems of this paper are as follows.

\medskip

{\bf Theorem \ref{qipairs}:}{\it
Let $\phi$ be a uniformly proper (bijective, by definition) map from $\JJ_1 \rightarrow \JJ_2$. There exists a quasi-isometry $q$ from $\Gamma_1$ to $\Gamma_2$ which pairs the sets $\JJ_1$ and $\JJ_2$ as $\phi$ does.}

The construction of the quasi-isometry $q$ proceeds by constructing a "coarse barycenter" of some infinite diameter sets (reminiscent of the celebrated measure-theoretic barycenter method discovered by Douady and Earle, and extended greatly by Besson, Courtois, Gallot
\cite{BCG} ).

We prove an analogous theorem
 for pairs $(X, \JJ)$ when $X$ is (strongly) hyperbolic relative to the collection $\JJ$. \\

{\bf Theorem \ref{qipairs-relhyp}:}{\it
Let $X_i$ be (strongly) hyperbolic relative to collections $\JJ_i$ ($i = 1, 2$).
Let $\phi$ be a uniformly proper (bijective, by definition) map from $\JJ_1 \rightarrow \JJ_2$. There exists a quasi-isometry $q$ from $X_1$ to $X_2$ which pairs the sets $\JJ_1$ and $\JJ_2$ as $\phi$ does.}

\smallskip

As a Corollary of Theorem \ref{qipairs-relhyp} and work of Hruska and Kleiner \cite{hruska-kleiner}, we deduce relative rigidity for
pairs $(X, \JJ)$ where $X$ is a CAT(0) space with isolated flats and  $\JJ$ is  the collection of maximal flats.

\smallskip

The third main theorem of this paper is an analog for higher
 rank symmetric spaces.\\

{\bf Theorem \ref{qipairs-ss}:}{\it
Let $X_i$ be symmetric spaces of non-positive curvature, and  $\JJ_i$ be equivariant collections of lifts of a maximal torus in a compact locally symmetric space modeled on $X_i$ ($i = 1, 2$).
Let $\phi$ be a uniformly proper (bijective, by definition) map from $\JJ_1 \rightarrow \JJ_2$. There exists a quasi-isometry $q$ from $X_1$ to $X_2$ which pairs the sets $\JJ_1$ and $\JJ_2$ as $\phi$ does.}

In fact, combining Theorem \ref{qipairs-ss} with the quasi-isometric rigidity theorem of Kleiner-Leeb \cite{kleiner-leeb} and Eskin-Farb \cite{eskin-farb-jams}, 
we may upgrade the quasi-isometry of Theorem \ref{qipairs-ss}  to an
isometry.

\medskip

Let $G_i, H_i$ ($i=1, 2$) be hyperbolic groups and quasiconvex subgroups respectively. In
Section 1.3, we shall construct simplicial complexes (termed $C$-complexes) from the incidence relations determined by
the cosets of $H_i$.
Let $C(G_i,H_i)$ be the $C$-complexes associated with the pairs $(G_i, H_i)$.   Roughly speaking, the vertices of $C(G_i,H_i)$ are the translates $g_i^j\Lambda_i$ of $\Lambda_i$ by distinct coset representatives $g_i^j$ and the $(n-1)$-cells are $n$-tuples $\{ g_1^1\Lambda , \cdots , g_1^n\Lambda \}$ of distinct translates such that $\cap_1^n g_1^i\Lambda \neq \emptyset$.

\medskip

{\bf Theorem \ref{isomccx}: }{\it
Let $\phi : \JJ_1 \rightarrow \JJ_2$ be a uniformly proper map. Then $\phi$ induces an isomorphism of $C(G_1,H_1)$ with $C(G_2,H_2)$.}

\medskip

On the way towards proving Theorems \ref{qipairs} and \ref{isomccx}, we prove two Propositions characterizing quasiconvexity. These might be of independent interest. The first is in terms of the Hausdorff topology on the collection $C_c^0(\partial G)$, which is the collection of closed subsets of $\bdy G$ having more than one point. 

\medskip

{\bf Proposition \ref{qc=discrete}:} {\it
Let $H$ be a subgroup of a hyperbolic group $G$ with limit set $\Lambda$. Let $\LL$ be the collection of translates of $\Lambda$ by elements of distinct cosets of $H$ (one for each coset). Then $H$ is quasiconvex if and only if $\LL$ is a discrete subset of $C_c^0(\partial G)$.}

\medskip

The second characterization is in terms of strong relative hyperbolicity.

\begin{defn} A subgroup $H$ of a group $G$ is said to be malnormal if for all $g \in G \setminus H$,
$gHg^{-1} \cap H$ is trivial. A subgroup $H$ of a group $G$ is said to be almost
malnormal if for all $g \in G \setminus H$,
$gHg^{-1} \cap H$ is finite. 
\end{defn}

It was pointed out to us
by the referee that the following result 
follows from work of Farb \cite{farb-relhyp},
Bowditch (\cite{bowditch-relhyp} Theorem 7.11) and 
Drutu-Sapir (\cite{drutu-sapir} Lemma 4.15). We shall 
include a  proof for completeness.

\medskip

{\bf Proposition \ref{rh2}:}\cite{farb-relhyp} 
\cite{bowditch-relhyp} \cite{drutu-sapir} {\it
Let $G$ be a hyperbolic group and $H$ a subgroup.
Then $G$ is strongly relatively hyperbolic with respect to $H$ if and only if 
 $H$ is  a malnormal quasiconvex subgroup.}

\smallskip

The prototypical example is that of (fundamental groups of) a closed hyperbolic manifold with a totally geodesic {\em embedded} submanifold.

\medskip

Finally, we give an intrinsic or dynamic reformulation of Theorems \ref{qipairs} and \ref{isomccx} following Bowditch \cite{bowditch-jams}, which makes use of the existence of a cross-ratio on the boundary of a hyperbolic group. The cross-ratio in turn induces a pseudometric on the collection $\LL$ of translates of $\Lambda$. 

\medskip

{\bf Theorem \ref{dyna}:}  {\it
Let $G_1, G_2$ be uniform convergence  (hence hyperbolic) groups acting on compacta $M_1, M_2$ respectively. Also, let $\AA_i$ (for $i = 1, 2$) be $G_i$-invariant annulus systems and let $(..|..)_i$ denote the corresponding annular cross-ratios.\\
 Let $H_1, H_2$ be subgroups of $G_1, G_2$ with limit sets $\Lambda_1, \Lambda_2$. Suppose that the set $\LL_i$ of translates of $\Lambda_i$ (for $i = 1, 2$) by essentially distinct elements of $H_i$ in $G_i$ forms a discrete subset  of $C_c^0(M_i)$. \\
Also assume that there exists a bijective function $\phi : \LL_1 \rightarrow \LL_2$ and  that this pairing is uniformly proper with respect to the cross-ratios $(..|..)_1$ and $(..|..)_2$.  Then \\
1) $H_i$ is quasiconvex in $G_i$ \\
2) There is a homeomorphism $q: M_1 \rightarrow M_2$ which pairs $\LL_1$ with $\LL_2$ as $\phi$ does. Further, $q$ is uniformly proper with respect to the cross-ratios $(..|..)_1$ and $(..|..)_2$ on $M_1$, $M_2$ respectively.\\
3) $q$ (and hence also $\phi$) induces an isomorphism of $C$-complexes $C(G_1,H_1)$ with $C(G_2,H_2)$.}

\medskip 

\noindent {\bf Acknowledgements:} My interest in relative hyperbolicity and quasi-isometric rigidity is largely due to  Benson Farb. It is a pleasure to acknowledge his help, support and camaraderie, both mathematical and personal. 
I would also like to thank the referee 
  for suggesting several corrections and for providing additional references.

\subsection{Relative Hyperbolicity and Electric Geometry}

We start off by fixing notions and notation.
Let $G$ (resp. $X$)
be a hyperbolic group (resp. a hyperbolic metric space)
with Cayley graph (resp. a net) $\Gamma$  equipped with a 
word-metric (resp. a simplicial metric) $d$. 

Here a net $\NN$ is a collection of distinct
points $x_i \in X$
such that there exist $0 < C_1 < C_2$ such that \\
1) $d(x_i, x_j) \geq C_1$ for all $i \neq j$\\
2) For all $x \in X$, there exists $x_i \in \NN$
such that $d(x_i, x) \leq C_2$\\

For the net $\NN$ we construct a graph $G_N$ with edges corresponding to pairs $x_i \neq x_j$ such that
$d(x_i,x_j) \leq 4C_2$.  The simplicial metric on $\NN$ is obtained
by declaring that each edge of $G_N$ has length one.

Let the 
{\bf Gromov boundary}
 of  $\Gamma$ be denoted by $\partial G$.
(cf.\cite{GhH}).

We shall have need for the fact that for hyperbolic metric spaces (in the sense of Gromov \cite{gromov-hypgps})
the notions of quasiconvexity and qi embeddings coincide  \cite{gromov-hypgps}.

We shall now recall certain notions of relative
hyperbolicity due to Gromov \cite{gromov-hypgps} and
 Farb \cite{farb-relhyp}.

Let $X$ be a path metric space. A collection of closed
 subsets $\HH = \{ H_\alpha\}$ of $X$ will be said to be {\bf uniformly
 separated} if there exists $\epsilon > 0$ such that
$d(H_1, H_2) \geq \epsilon$ for all distinct $H_1, H_2 \in \HH$.

The {\bf electric space} (or coned-off space) $\hhat{X}$
corresponding to the
pair $(X,\HH )$ is a metric space which consists of $X$ and a
collection of vertices $v_\alpha$ (one for each $H_\alpha \in \HH$)
such that each point of $H_\alpha$ is joined to (coned off at)
$v_\alpha$ by an edge of length $\half$.

\begin{defn} \cite{farb-relhyp} \cite{bowditch-relhyp}
Let $X$ be a geodesic metric space and $\HH$ be a collection of
 uniformly separated subsets. Then $X$ is said to be
{\bf weakly hyperbolic} relative to the collection $\HH$, if the
electric space $\hhat{X}$ is hyperbolic.
\end{defn}

\begin{lemma} (See Bowditch \cite{bowditch-relhyp},
generalizing Lemma 4.5 and Proposition 4.6 of Farb
  \cite{farb-relhyp})
Given $\delta , C, D$ there exists $\Delta$ such that
if $X$ is a $\delta$-hyperbolic metric space with a collection
$\mathcal{H}$ of $C$-quasiconvex $D$-separated sets.
then, the electric space
  $\hhat{X}$  is $\Delta$-hyperbolic, i.e. $X$ is weakly hyperbolic relative to the collection $\HH$. \\
\label{farb1A}
\end{lemma}

{\bf Definitions:} Given a collection $\mathcal{H}$
of $C$-quasiconvex, $D$-separated sets and a number $\epsilon$ we
shall say that a geodesic (resp. quasigeodesic) $\gamma$ is a geodesic
(resp. quasigeodesic) {\bf without backtracking} with respect to
$\epsilon$ neighborhoods if $\gamma$ does not return to $N_\epsilon
(H)$ after leaving it, for any $H \in \mathcal{H}$. 
A geodesic (resp. quasigeodesic) $\gamma$ is a geodesic
(resp. quasigeodesic) {\bf without backtracking} if it is a geodesic
(resp. quasigeodesic) without backtracking with respect to
$\epsilon$ neighborhoods for some $\epsilon \geq 0$.

Electric $P$-quasigeodesics without backtracking
are said to have similar intersection patterns if for
 $\beta , \gamma$ 
   electric $P$-quasigeodesics without backtracking
  both joining $x, y$, the following are satisfied.
\begin{enumerate}
\item {\it Similar Intersection Patterns 1:}  if
  precisely one of $\{ \beta , \gamma \}$ meets an
  $\epsilon$-neighborhood $N_\epsilon (H_1)$
of an electrocuted quasiconvex set
  $H_1 \in \mathcal{H}$, then the length (measured in the intrinsic path-metric
  on  $N_\epsilon (H_1)$ ) from the entry point
  to the 
  exit point is at most $D$. \\
\item {\it Similar Intersection Patterns 2:}  if
 both $\{ \beta , \gamma \}$ meet some  $N_\epsilon (H_1)$
 then the length (measured in the intrinsic path-metric
  on  $N_\epsilon (H_1)$ ) from the entry point of
 $\beta$ to that of $\gamma$ is at most $D$; similarly for exit points. \\
\end{enumerate}

\begin{defn} \cite{farb-relhyp} \cite{bowditch-relhyp}
Let $X$ be a geodesic metric space and $\HH$ be a collection of
mutually disjoint uniformly separated subsets such that
 $X$ is 
{\em weakly hyperbolic} relative to the collection $\HH$. If any pair of $P$-
electric quasigeodesics without backtracking starting and ending at
the same point have {\em similar intersection patterns} with {\em
  horosphere-like sets} (elements of $\HH$) then quasigeodesics
are said to satisfy {\bf Bounded Penetration} and 
 $X$ is said to be
{\bf strongly hyperbolic} relative to the collection $\HH$.
\end{defn}

\begin{defn} \cite{brahma-ibdd} A collection $\mathcal{H}$ of uniformly
$C$-quasiconvex sets in a $\delta$-hyperbolic metric space $X$
is said to be {\bf mutually D-cobounded} if 
 for all $H_i, H_j \in \mathcal{H}$, $\pi_i
(H_j)$ has diameter less than $D$, where $\pi_i$ denotes a nearest
point projection of $X$ onto $H_i$. A collection is {\bf mutually
  cobounded } if it is mutually D-cobounded for some $D$. 
\end{defn}

{\em Mutual coboundedness} was proven by Farb
for horoballs  in finite volume Hadamard manifolds of pinched
negative curvature in Lemma 4.7 of
\cite{farb-relhyp}. The following generalization
is due to Bowditch \cite{bowditch-relhyp}.

\begin{lemma}  (See
Bowditch \cite{bowditch-relhyp} Lemma 7.13 for a proof)
Suppose $X$ is a $\delta$-hyperbolic metric space with a collection
$\mathcal{H}$ of $C$-quasiconvex $K$-separated $D$-mutually cobounded
subsets. Then $X$ is strongly hyperbolic relative to the collection $\HH$. 
\label{farb2A}
\end{lemma}

\medskip

Gromov gave a different definition of {\bf strong relative hyperbolicity}. 
We give a condition below that is equivalent to a special
case of Gromov's definition.
Let $X$ be a geodesic metric space with a collection $\mathcal{H}$ of uniformly separated subsets $\{ H_i \}$. The hyperbolic cone $cH_i$ is the product of $H_i$ and the non-negative reals
 $H_i \times \RR_+$, equipped with the metric of the type $2^{-t}ds^2 + dt^2$. More precisely, $H_i \times \{ n \}$ is given the path metric of $H_i$ scaled by $2^{-n}$. The $\RR_+$ direction is given the standard Euclidean metric.
Let $X^h$ denote $X$ with hyperbolic cones $cH_i$ glued to it along $H_i$'s. $X^h$ will be referred to as the {\em hyperbolically coned off $X$}. This is to be contrasted with the coned off space $\hat{X}$ in Farb's definition.

\begin{defn}
$X$ is said to be strongly hyperbolic relative to the collection $\mathcal{H}$ in the sense of Gromov if the 
hyperbolically coned off space $X^h$ is a hyperbolic metric space.
\label{gromovrh}
\end{defn}

The equivalence of the two notions of strong relative hyperbolicity was proven by Bowditch in \cite{bowditch-relhyp}.

\begin{theorem} ( Bowditch \cite{bowditch-relhyp} )
$X$ is  strongly hyperbolic relative to a collection $\mathcal{H}$ of uniformly separated subsets $\{ H_i \}$ in the sense of Gromov if and only if $X$ is  strongly hyperbolic relative to the collection $\mathcal{H}$  in the sense of Farb. 
\label{farb=gromov}
\end{theorem}

\subsection{Height of Subgroups and C-Complexes}

The notion of height of a subgroup was introduced by Gitik, Mitra, Rips and Sageev in \cite{GMRS} and further developed by the author in \cite{mitra-ht}.

\medskip

\begin{defn}
Let $H$ be a subgroup of a  group $G$. We say that
the elements $\{g_i |1 \le i \le n\}$ of $G$ are 
essentially distinct if $Hg_i \neq Hg_j$ for $i \neq j$.
Conjugates of $H$ by essentially distinct elements are called
 essentially distinct conjugates. 
\end{defn}

\medskip

Note that we are abusing notation slightly here, as a conjugate of $H$
by
an element belonging to  the 
 normalizer of $H$ but not belonging to
 $H$ is still essentially distinct from $H$.
Thus in this context a conjugate of $H$ records (implicitly) the conjugating
element.

\medskip

\begin{defn}
We say that the height of an infinite subgroup $H$ in $G$ is $n$  if  
there exists a collection of $n$ essentially distinct conjugates
of $H$ such that the intersection of all the elements of the collection is 
infinite  and $n$ is maximal possible. We define the height of a finite 
subgroup  to be $0$. We say that the width of an infinite subgroup $H$ in $G$ is $n$  if  
there exists a collection of $n$ essentially distinct conjugates
of $H$ such that the {\bf pairwise} intersection of  the elements of the collection is 
infinite  and $n$ is maximal possible.
\end{defn}

The  main  theorem of \cite{GMRS} states:

\begin{theorem}
If $H$ is a quasiconvex subgroup of a hyperbolic group $G$,then
$H$ has finite height and finite width.
\label{gmrs}
\end{theorem}

In this context, a theorem we shall be needing several times is the following result from \cite{GMRS}
that is proved using a result of  Short \cite{short}.

\begin{theorem}  (Lemma 2.6 of \cite{GMRS})
Let $G$ be  a hyperbolic group and $H_i$ (for $i = 1 \cdots k$ ) be quasiconvex subgroups with limit sets $\Lambda_i$, 
$i = 1 \cdots k$. Then $\cap H_i$ is a quasiconvex subgroup with limit set  $\cap \Lambda_i$.
\label{shortthm}
\end{theorem}

We now proceed to define a simplicial
complex $C(G,H)$ for a  group $G$ and $H$ a  subgroup. For $G$ hyperbolic and $H$ quasiconvex, we give below three equivalent descriptions of a complex $C(G,H)$. In this case, let $\partial G$ denote the boundary of $G$, $\Lambda$ the limit set of $H$, and $J$ the join of $\Lambda$. \\
1) Vertices ( $0$-cells ) are  conjugates of $H$ by essentially distinct elements, and $(n-1)$-cells are $n$-tuples $\{ g_1H, \cdots , g_nH \}$ of distinct cosets such that $\cap_1^n g_iHg_i^{-1}$ is infinite (in fact by Theorem \ref{shortthm} an infinite quasiconvex subgroup of $G$).\\
2) Vertices ( $0$-cells ) are translates of $\Lambda$ by essentially distinct elements, and $(n-1)$-cells are $n$-tuples $\{ g_1\Lambda , \cdots , g_n\Lambda \}$ of distinct translates such that $\cap_1^n g_i\Lambda \neq \emptyset$.\\
3) Vertices ( $0$-cells ) are  translates of $J$ by essentially distinct elements, and $(n-1)$-cells are $n$-tuples $\{ g_1J, \cdots ,g_nJ \}$ of distinct translates such that $\cap_1^n g_iJ$ is infinite. 

\medskip

We shall refer to the complex $C(G,H)$ as the {\bf C-complex}  for the pair $G, H$. ({\bf C} stands for ``coarse" or ``\v{C}ech" or ``cover", since $C(G,H)$ is like a coarse nerve of a cover, reminiscent of constructions in Cech cochains.) {\em Note that if $h(H)$ denote the height of $H$, then $(h(H)+1)$ is the dimension of the $C$-complex $C(G,H)$. Also, if $w(H)$ denote the width of $H$, then $w(H) = w$ is equal to the size of the largest complete graph $K_w$ that is embeddable in $C(G,H)$.} If $C(G,H)$ is connected then its one-skeleton is closely related to the coned off space $\hat{\Gamma }$ with an appropriately chosen set of generators. 

This definition is inspired by that of the Tits complex for a non-uniform lattice in a higher rank symmetric space. Related constructs in the context of codimension 1 subgroups also occur in work of Sageev \cite{sageev-th} where he constructs cubings.

\section{Characterizations of Quasiconvexity}

Let $G$ be a hyperbolic group.
Let $C_c(\partial G)$ denote the collection of closed subsets of the boundary $\partial G$ equipped with the Hausdorff topology. Let $C_c^0(\partial G) \subset C_c(\partial G)$  denote the subset obtained from $C_c(\partial G)$ by removing the singleton sets $\{ \{ x \}: x \in \partial G \}$. Next fix a subgroup $H \subset G$ with limit set $\Lambda \subset \partial G$. Consider the $G$-invariant collection
$\mathcal{L} = $$\{$ $ g \Lambda $ $\}$
$ \subset C_c^0(\partial G)$ with $g$ ranging over {\em distinct cosets} (one for each coset) of $H$ in $G$. Note that $\LL$ is (strictly speaking) a {\em multi-set} as distinct elements of $\LL$ may denote the same element of $C_c^0(\partial G)$ in case two distinct translates of $\Lambda$ coincide. One extreme case is when $\Lambda = \partial G$, though $H$ is of infinite index in $G$ (e.g. if $H$ is normal of infinite index in in $G$.) Then $\LL$ consists of infinitely many copies of $\Lambda$. 

\begin{defn}
The {\bf join} $J( \Lambda )$ of $\Lambda$ is defined as the union of all bi-infinite geodesics whose end-points lie in $\Lambda$
\end{defn}

It is easy to see that $J(\Lambda )$ is $2 \delta$-quasiconvex if $G$ is $\delta$-hyperbolic. In fact this is true for any subset $\Lambda$ of the boundary of a $\delta$-hyperbolic metric space $X$ (no equivariance is necessary). For $\Lambda$ the limit set of $H$, $J(\Lambda )$ is $H$-invariant. The {\em visual diameter}  $dia_{\partial G}(\Lambda ) $ of a subset $\Lambda$ of $\partial G$ is the same as the diameter in the metric on $\partial G$ obtained from the Gromov inner product. (See \cite{GhH} Chapter
7 for details about the visual metric on $\partial G$.)

\subsection{Limit Sets and Quasiconvexity}

The next Lemma follows directly from the fact that sets with  visual diameter bounded below contain points with Gromov inner product bounded above  and conversely\cite{GhH}.

\begin{lemma}
 For all $\epsilon > 0$ there exists $N$ such that if the diameter $dia_{\partial G}(\Lambda ) \geq \epsilon$ for a    closed subset $\Lambda$ of $\partial G$, then there exists $p\in J(\Lambda )$ such that $d(p,1) \leq N$. Conversely, 
for all $N > 0$ there exists $\ep > 0$ such that if
there exists $p\in J(\Lambda )$ with $d(p,1) \leq N$, then
$dia_{\partial G}(\Lambda ) \geq \epsilon$.
\label{vdia}
\end{lemma}

The next Proposition gives our first characterisation of quasiconvex subgroups of a hyperbolic group.

\begin{prop} {\bf (Characterization of Quasiconvexity I)}
Let $H$ be a subgroup of a hyperbolic group $G$ with limit set $\Lambda$. Let $\LL$ be the collection of translates of $\Lambda$ (counted with multiplicity)
 by elements of distinct cosets of $H$ (one for each coset). Then $H$ is quasiconvex if and only if $\LL$ is a discrete subset of $C_c^0(\partial G)$.
\label{qc=discrete}
\end{prop}

\noindent {\bf Proof:} Suppose $H$ is quasiconvex. We want to show that $\LL$ is a discrete subset of $C_c^0(\partial G)$. Thus it suffices to show that any limit of elements of $\LL$ is a singleton set. This in turn follows from the following.

\medskip

\noindent {\bf Claim:} For all $\epsilon > 0$, $\LL_\epsilon = \{ L_i \in \LL : dia_{\partial G} (L_i) \geq \ep \}$ is finite. \\
{\bf Proof of Claim:} Let $N = N(\ep )$ be as in Lemma \ref{vdia}. Since $ dia_{\partial G} (L_i) \geq \ep $, therefore by Lemma \ref{vdia}, there exists $p_i \in J(L_i)$ such that $d_G(p_i,1) \leq N$. Also,  there exists  $K > 0$ depending on $\delta$ (recall that $J(L_i)$ is $2 \delta$-qc) and the quasiconvexity constant of $H$ such that if  $L_i = g_i \Lambda$, then there exists $h_i \in H$ with $d_G (p_i, g_ih_i) \leq K$. Hence, $d_G(1, g_ih_i) \leq K+N$. Since $G$ is finitely generated, the number of such elements $g_ih_i$ is finite. Since $g_i$ are picked from distinct cosets of $H$, we conclude that the set $\LL_\epsilon$ is finite. $\Box$

Conversely, suppose that $H$ is not quasiconvex. 
Assume, without loss of generality, that a finite 
generating set of $H$ is contained in a 
finite 
generating set of $G$ and that $\Gamma_H , \Gamma_G$
are Cayley graphs with respect to these generating sets.
Then there exist $p_i \in J(\Lambda )$ such that $d_G(p_i, \Gamma_H) \geq i$. Translating by an appropriate element of $H$, we may assume that $d_G(p_i, \Gamma_H) = d_G(p_i, 1) \geq i$. Further, we may assume (by passing to a subsequence if necessary) that the sequence $d_G(p_i, 1)$ is monotonically increasing. Then $p_i^{-1}J(\Lambda )$ has limit set $p_i^{-1} \Lambda$. Further, as $p_i \in J(\Lambda )$, therefore, 
$1 \in p_i^{-1}J(\Lambda )$.
Since $J(\Lambda )$ is $2\delta$-qc, so is $p_i^{-1}J(\Lambda )$ for all $i$. Hence, there exists $\ep > 0$ by Lemma \ref{vdia} such that $dia_{\partial G} p_i^{-1}J(\Lambda ) \geq \ep$. Since $d_G(p_i, 1)$ is monotonically strictly increasing, we conclude that $p_i$'s lie in distinbct cosets of $H$. Further, since $C_c(\partial G)$ is compact, we conclude that the collection $p_i^{-1}J(\Lambda )$ has a convergent subsequence, converging to a subset of diameter greater than or equal to $\ep$. Therefore, the collection $\LL$ is {\em not} a discrete subset (strictly speaking a multiset) of $C_c^0(\partial G)$. $\Box$

We next prove a result about projections of $J(L_i)$ on 
$J(L_j)$. 
We start off with an elementary fact about hyperbolic metric spaces. See \cite{mitra-trees} for a proof.

\begin{lemma} \cite{mitra-trees}
Given $\delta > 0$,
 there exist $D, C_1, k, \ep$ such that  if $a, b, c, d$
are points of a $\delta$-hyperbolic metric space $(Z,d)$,
with ${d}(a,[b,c])={d}(a,b)$,
${d}(d,[b,c])={d}(c,d)$ and ${d}(b,c)\geq{D}$
then $[a,b]\cup{[b,c]}\cup{[c,d]}$ lies in a $C_1$-neighborhood of
any geodesic joining 
$a, d$ and is a $(k, \ep)$-quasigeodesic.
\label{perps}
\end{lemma}

Assume that $H$ is quasiconvex and that $L_k$ is the limit set $g_k\Lambda$ of $g_kH$. Let $P_{j}$ denote the nearest point projection of $\Gamma_G$ onto $J(L_j)$. Also, let $H_k = g_k \Gamma_H $ be the left translate of $\Gamma_H$ by $g_k$.

\begin{prop}
There exists $K >0$ such that $P_j(\Gamma_{H_i})$ lies in a $K$-neighborhood of $J(L_i \cap L_j)$ if $(L_i \cap L_j) \neq \emptyset$. Else, $P_j(\Gamma_{H_i})$ has diameter less than $K$. 
\label{cobdd}
\end{prop}

\noindent {\bf Proof:} Since $J(L_i)$ is $2\delta$-qc and $H$ is quasiconvex, it suffices to show that $P_j(J(L_i))$ lies in a $K$-neighborhood of $J(L_i \cap L_j)$
if the latter is non-empty. By $G$-equivariance, we may assume that $L_j = \Lambda$ and 
$g_i = 1$. We represent $P_j$ by $P$ in this case.

First note that by  Theorem \ref{shortthm}, 
$H_i \cap H_j$ is quasiconvex and the limit set of 
$H_i \cap H_j$ is $L_i \cap L_j$. Also, 
$J(L_i \cap L_j) \subset J(L_i)$. 

Let $a, b \in J(L_i)$. Let $P(a) = c, P(b) = d$. Let $D, C_1, k, \ep$ be as in Lemma \ref{perps}. If $d_G(c,d) \geq D$, then $[a,c]\cup [c,d]\cup [d,b]$ is a $(k, \ep )$-quasigeodesic lying in a $C_1$ neighborhood of $[a,b]$. Since $J(L_i)$, $J(L_j)$ are both $2\delta$-qc, $[a,b]$ lies in a $2\delta$ -neighborhood of 
$J(L_i)$, and $[c,d]$ lies in a $2\delta$ -neighborhood of 
$J(L_j)$. In particular $c,d$ lie in  a $(C_1 + 2 \delta )$-neighborhood of $J(L_j)$. Translating by an element of $H$, we may assume that $c = 1$. (Note that the argument in this
paragraph works independent of whether
$J(L_i) \cap J(L_j)$ is empty or not.

We proceed now by contradiction. Suppose there exists a sequence of $L_i$'s and $b_i \in J(L_i)$ such that $P (b_i) = d_i$ lies at a distance greater than $i$ from 
$J(L_i \cap L_j) $ (resp. $c=1$) according as 
$J(L_i) \cap J(L_j)$ is non-empty or
empty. This shows that the sequence $L_i$ has a limit point on $\Lambda$ disjoint from $L_i \cap \Lambda$ for all $i$ and further that $J(L_i)$ passes through a bounded neighborhood of $1$. Hence the sequence $L_i$  is not discrete in $C_c^0(\partial G)$. This contradicts Proposition \ref{qc=discrete} and proves our claim. $\Box$

\subsection{Quasiconvexity and Relative Hyperbolicity}

As an immediate corollary of Proposition \ref{cobdd} in conjunction with Theorem \ref{shortthm} of Short \cite{short}, we immediately conclude

\begin{cor}
Let $H$ be a malnormal quasiconvex subgroup of a hyperbolic group $G$ with Cayley graph $\Gamma$ and limit set $L$. Then the set of joins $\JJ$ of distinct translates of $L$ is a uniformly cobounded collection of uniformly quasiconvex sets in $\Gamma$.
\label{cobddcor}
\end{cor}

Combining Lemma \ref{farb2A} with Corollary \ref{cobddcor} above, we have the following
Proposition due to Bowditch \cite{bowditch-relhyp}.

\begin{prop}{\bf (Characterization of Quasiconvexity II)
\cite{bowditch-relhyp}}
Let $H$ be a malnormal quasiconvex subgroup of a hyperbolic group $G$. Then $G$ is strongly relatively hyperbolic with respect to $H$.
\label{rh1}
\end{prop}

In fact the converse to Proposition \ref{rh1} is also true.
We came to learn from the referee that this follows by combining work of Farb \cite{farb-relhyp}, Bowditch \cite{bowditch-relhyp} and Drutu-Sapir \cite{drutu-sapir}. We provide
a proof below for completeness (and because it is easily
done).

Malnormality of strongly relatively hyperbolic subgroups is due to Farb \cite{farb-relhyp}. In fact this does not require $G$ to be hyperbolic.

\begin{lemma} (Farb \cite{farb-relhyp}) 
Let $G$ be strongly relatively hyperbolic with respect to $H$. Then $H$ is malnormal in $G$.
\label{mn}
\end{lemma}

It remains to show that $H$ is quasiconvex if  a hyperbolic group $G$ be strongly relatively hyperbolic with respect to $H$. We use Gromov's definition of strong relative hyperbolicity. Attach hyperbolic cones $cH$ to distinct translates of $\Gamma_H$ in $\Gamma_G$ to obtain the hyperbolically coned off Cayley graph $\Gamma_G^h$. Then 
$\Gamma_G^h$ is hyperbolic by Gromov's definition. 

If $H$ is not quasi-isometrically embedded in $G$ then for all $i \in \natls$, there exist $p_{i1}, p_{i2} \in \Gamma_H$ such that 
$$d_H(p_{i1}, p_{i2}) \geq id_G(p_{i1}, p_{i2})$$.
Also from the metric $d_{cH}$ on $cH$, we find that $d_{cH}(p_{i1}, p_{i2})$ is of the order of $log_2d_H(p_{i1}, p_{i2}) $. Hence, we can further assume that 
$$d_H(p_{i1}, p_{i2}) \geq id_{cH}(p_{i1}, p_{i2})$$.
Join $p_{i1}, p_{i2}$ by shortest paths $\alpha_i, \beta_i$ in $cH$, $\Gamma_G$ respectively. Then $\alpha_i \cup \beta_i = \sigma_i$ is a closed loop in $\Gamma_G^h$ with total length $l(\sigma_i) = 
(d_{cH}(p_{i1}, p_{i2}) + d_G(p_{i1}, p_{i2}))$.
Therefore $il(\sigma_i) \leq 2d_H(p_{i1}, p_{i2}) $.

Since any (combinatorial) disk $D_i$ in $\Gamma_G^h$ spanning $\sigma_i$  must contain a path $\gamma_i$ in $\Gamma_H$ joining $p_{1i},p_{2i}$, therefore the area $A(D_i)$ of $D_i$ must be at least that of $N_1(\gamma_i)$, the $1$-neighborhood of $\gamma_i$ in $D_i$.

Therefore there exists $C > 0$ such that for all $i$, 
$$A(D_i) \geq A(N_1(\gamma_i)) \geq \frac{d_H(p_{i1}, p_{i2})}{C} \geq \frac{il(\sigma_i)}{2C} $$

Since $i$ is arbitrary, this shows that $\Gamma_G^h$ cannot satisfy a linear isoperimetric inequality. Hence
$\Gamma_G^h$ cannot be a hyperbolic metric space. This is a contradiction. Hence $H$ must be quasi-isometrically embedded in $G$. Hence (see for instance \cite{gromov-hypgps} ), $H$ is quasiconvex in $G$. This completes our proof of the following characterisation of strongly relatively hyperbolic subgroups of hyperbolic groups.

\begin{prop}
Let $G$ be a hyperbolic group and $H$ a subgroup.
Then $G$ is strongly relatively hyperbolic with respect to $H$ if and only if 
 $H$ is  a malnormal quasiconvex subgroup.
\label{rh2}
\end{prop}

\section{Relative Rigidity}

\subsection{Pairing of Limit Sets by Quasi-isometries}

We now consider two hyperbolic groups $G_1, G_2$ with quasiconvex subgroups $H_1, H_2$, Cayley graphs $\Gamma_1, \Gamma_2$. Let $\LL_j$  for $j = 1, 2$ denote the collection of translates of limit sets (counted with multiplicity as before) of $H_1, H_2$ in $\partial G_1, \partial G_2$ respectively. Individual members of  the collection $\LL_j$ will be denoted as $L^j_i$.
Let $\JJ_j$ denote the collection 
$\{ J_i^j = J(L_i^j): L_i^j \in \LL_j \}$.
Following Schwarz \cite{schwarz-inv}, we define:

\begin{defn} A bijective map $\phi$ from $\JJ_1 \rightarrow \JJ_2$  is said to be uniformly proper  if there exists a function $f: \natls \rightarrow \natls$ such that \\
1) $d_{G_1} (J(L_i^1) ,J(L_j^1)) \leq n \Rightarrow d_{G_2} (\phi(J(L_i^1)) ,\phi(J(L_j^1))) \leq f(n)$ \\
2)  $ d_{G_2} (\phi(J(L_i^1)) ,\phi(J(L_j^1)))\leq n \Rightarrow d_{G_1} (J(L_i^1) ,J(L_j^1)) \leq f(n)$. 

When $\JJ_i$ consists 
of all singleton subsets of $\Gamma_1, \Gamma_2$, we shall
refer to $\phi$ as a uniformly proper map from 
$\Gamma_1$ to $ \Gamma_2$.
\label{uproper}
\end{defn}

\noindent {\bf Note:} We observe that if $\JJ_i$ is just the collection of singleton sets in $\Gamma_i$, then a uniformly proper map between $\JJ$'s is 
 the same as a quasi-isometry between $\Gamma_i$'s. This can be seen by putting $n = 1$ in conditions 1 and 2 above
and then using the fact that graphs have edge length one.
Hence what is important here is that $\JJ$'s are infinite diameter sets.  

\begin{defn}
A map $q$ from $\Gamma_1$ to $\Gamma_2$ is said to pair the sets $\JJ_1$ and $\JJ_2$ as $\phi$ does if there exists a function $h: \natls \rightarrow \natls$ such that \\
$d_G (p ,J^1_j)) \leq n \Rightarrow d_G (q(p) ,
\phi (J(L_j^1))) \leq h(n)$.
\end{defn}

The following Lemma generalises Lemma 3.1 of Schwarz \cite{schwarz-inv}, where the result is proven in the special case of a symmetric pattern of geodesics in $\Hyp^n$. The referee
pointed out to us that the Lemma follows from Lemma 7
of \cite{niblo-reeves} by Niblo and Reeves.

\begin{lemma}
For $M, m > 0$, there exists $R > 0$, such that the following holds. \\
Let $L_1, \cdots , L_M$ be distinct translates of the limit set of a quasiconvex subgroup $H$ of a hyperbolic group $G$, such that $d_G(J_i,J_j) \leq m$ for all $i, j = 1 \cdots , M$ and $J_i = J(L_i)$. Then there exists a ball of radius $R$ meeting $J_i$ for all $i = 1 \cdots , M$.
\label{3.1}
\end{lemma}

\noindent {\bf Proof:} If $\cap_1^M L_i \neq \emptyset$, choose any point $p \in J(\cap_1^M L_i )$. Then $B_1(p)$ intersects all $J_i$ and we are through.

Suppose therefore that $\cap_1^m L_i = \emptyset$. We proceed by induction on $M$. There exists $R_{M-1}$ such that a ball of radius $R_{M-1}$ meets $J_i$ for $i = 1 \cdots M-1$. 

We now proceed by contradiction. If no such $R$ exists for $M$, we have collections $\{ L_1^k, \cdots , L_M^k \}, k \in \natls$ such that a ball of radius $R_{M-1}$ meets $J_i^k, i = 1 \cdots M-1$ but no ball of radius $k$ meets $J_i^k, i = 1 \cdots M$. In particular, (since $J(\cap_1^{M-1} L_i^k) \subset \cap_1^{M-1} J_i^k)$), if $\cap_1^{M-1} L_i^k \neq \emptyset$, then $N_k(J(\cap_1^{M-1} L_i^k)) \cap J_i^M = \emptyset $.

For all $i, j, k$, choose points $p_{ij}^k \in J_i^k$ such that $d_G(p_{ij}^k, p_{ji}^k) \leq m$.

Assume by $G$-invariance of $\JJ$ that the ball of radius $R_{M-1}$ centered at $1 \in \Gamma_G$ meets $J_i^k , i = 1 \cdots M-1$. Therefore $J^k_M$ lies outside a $k$-ball about $1$. 

\medskip

Since the collection of $J_i$'s through $1$ is finite, therefore assume after passing to a subsequence if necessary, that \\
1) $\{ J^k_i \}_k$ is a constant sequence for $i = 1 \cdots M-1$. Hence, $\{ L^k_i \}_k$ is a constant sequence $L_i$ (say) for $i = 1 \cdots M-1$. \\
2) $p_{iM}^k \rightarrow p_{iM} \in \partial G$ for $i = 1 \cdots M-1$. Hence 
$p_{Mi}^k \rightarrow p_{iM} \in \partial G$. Further, by (1) above, $p_{iM} \in L_i$. \\
3) $L_M^k$ converges to a closed set $Z \subset \partial G$. By Proposition \ref{qc=discrete}, $Z$ must be a singleton set $\{ z \}$. \\
4) $J^k_M$ lies outside $B_k(1) \cup N_k(J(\cap_1^{M-1} L_i^k))$. If $\cap_1^{M-1} L_i \neq \emptyset$, then assume further by $G$-invariance, that $1 \in J(\cap_1^{M-1} L_i^k)$. Also, using Theorem \ref{shortthm} due to Short \cite{short}, and translating by an appropriate element of $\cap_1^{M-1} H_i^k$, we may assume that $1 \in J(\cap_1^{M-1} L_i^k)$ is closest to $J^k_M$. \\

Now, $p_{Mi}^k \in J^k_M$ and hence by (3) above, 
$p_{Mi}^k \rightarrow z \in \partial G$. Combining this with (2) above, we get $z = p_{iM}$ for all $i = 1 \cdots M-1$. 
Therefore, $z \in \cap_1^{M-1} L_i \neq \emptyset$.

But $d_G(1, J^k_M) = d_G(J(\cap_1^{M-1} L_i^k), J^k_M) \geq k$. 
Let $z_k \in J^k_M$ such that $d_G(1, J^k_M) = d_G(1,z_k) = d_G(J(\cap_1^{M-1} L_i^k), J^k_M) \geq k$. 

Then the Gromov inner product $(z_k,p_{iM}^k)_1$ is uniformly bounded above. Therefore $(z_k,p_{iM})_1$ is uniformly bounded above. Hence finally $(z,p_{iM})_1$ is bounded above. In particular $z \neq p_{iM}$. This is the contradiction that proves the Lemma. $\Box$

\smallskip

\noindent {\bf Definition of $q$} \\

\smallskip

Let $\phi$ be a uniformly proper (bijective, by definition) map from $\JJ_1 \rightarrow \JJ_2$. 
We shall now show
tha there exists a quasi-isometry $q$ from $\Gamma_1$ to $\Gamma_2$ which pairs the sets $\JJ_1$ and $\JJ_2$ as $\phi$ does.

 We will
 define a map $q: \Gamma_1 \rightarrow \Gamma_2$ which pairs $\JJ_1$ with $\JJ_2$ as $\phi$ does and prove that
$q$ is a quasi-isometry as promised.

Choose $K > 0$ such that the $K$ neighborhood $B_K(g)$ of $g \in \Gamma_1$
has greater than $w_2$ ( the  width of $H_2$ in $G_2$ )
$J^1_i$'s passing through it. 

Let $\JJ_{K,g}^j$ (for $j = 1, 2$ ) denote the collection of $J^j_i$'s passing through $N_K(g)$ for $g \in \Gamma_j, j = 1, 2$.

By the proof of Proposition \ref{qc=discrete}, there exists $M=M(K)$ (independent of $g \in \Gamma_1$)
such that at most $M$ $J^1_i $'s 
in $\JJ_{K,g}^1$ pass through 
$N_K(g)$. Since $\phi$ is a bijective pairing, $\phi ( \JJ_{K,g}^j)$ has at least $(w_2 + 1)$ and at most $M(K)$ elements in it.
By definition of $w_2$, and by  
Theorem \ref{shortthm} at least two of the limit sets of the $\phi (J^1_i)$'s are disjoint. Let $L_1^2$ and $L_2^2$ denote these limit sets. Hence, by Corollary \ref{cobddcor}, for any $K_1 \geq f(K)$, there exists $D$ such that the collection of points $$\{ p \in \Gamma_2 : d_2(p,J^1_2) \leq K_1, d_2(p,J^2_2) \leq K_1 \}$$ has diameter less than $D$.

Also, by uniform properness of $\phi$, $$d_2(\phi (J_m^1), \phi (J_n^1)) \leq f(2K)$$ for $J_m^1, J_n^1 $ passing through $N_K(g)$ (independent of $g$).

Summarising, \\
1) $L_1^2$ and $L_2^2$ are disjoint. \\
2) But, by Lemma \ref{3.1}, using $m = f(2K) $ and $M=M(K)$, there exists $R = R(K)$ and a ball of radius $R$ meeting each $\phi (J_i^1)$. \\
3) For any $K_1$, there exists $D $, such that $\{ p \in \Gamma_2 : d_2(p,J^1_2) \leq K_1, d_2(p,J^2_2) \leq K_1 \}$ has diameter less than $D$.In particular, we may choose $K_1 = R$. \\

\medskip

Define $q(g)$ to be the center of the ball of radius $R$ obtained in (2) above. By (3), $q(g)$ is thus defined upto a uniformly bounded amount of discrepancy for all $g \in \Gamma_1$. 

\begin{lemma}
$q$ is uniformly proper with respect to the metrics $d_1$, $d_2$.
\label{3.2}
\end{lemma}

\noindent {\bf Proof:} The proof is an almost exact replica of Lemma 3.2 of Schwarz \cite{schwarz-inv} and we content ourselves with reproducing the heuristics of his argument here. 

If $x, y$ are close in $\Gamma_1$, then the pairwise minimal distances between elements of $J_{K1}^x$ and $J_{K1}^y$ is uniformly bounded above. Hence, by Lemma \ref{3.1}, there exists a uniform upper bound to the radius of a minimal radius ball intersecting all elements of $\phi ( J_{K1}^x )$ as well as $\phi ( J_{K1}^y )$. Also, since the center $w$ of such a ball is defined upto a bounded amount of discrepancy, it must be at a bounded distance from both $q(x)$ as well as $q(y)$. Hence $d_2(q(x),q(y))$ is uniformly bounded, i.e. close.

Conversely, suppose that $q(x), q(y)$ are close. First, by Lemma \ref{3.1}, there exists a uniform upper bound $R$ on radius of minimal radius balls $B_1, B_2$ centered at $q(x), q(y)$, intersecting all elements of $\phi (\JJ^x_{K1}), \phi (\JJ^y_{K1})$ respectively. Then the $(R + d_2(q(x),q(y)))$ ball about $q(x)$ (or $q(y)$) meets every element of 
$\phi (\JJ^x_{K1})$ as well as $ \phi (\JJ^y_{K1})$. Since $\phi$ is uniformly proper, this means that there is a uniform upper bound on the minimal radius of a ball meeting every element of 
$ (\JJ^x_{K1})$ as well as $  (\JJ^y_{K1})$. As before, $d_1(x,y)$ is uniformly bounded, i.e. $x, y$ are close. $\Box$

\medskip

Similarly, we can construct $q^{-1}$ using the bijective pairing $\phi^{-1}$ such that $q^{-1}$ is uniformly proper.
Also, from Lemma \ref{3.1}  $q, q^{-1}$ composed with each other in either direction is close to the identity.

Since $\phi$ pairs $\LL_1$, $\LL_2$ bijectively and is uniformly proper from $\JJ_1$ to $\JJ_2$, therefore by invariance of $\JJ_2$ under $G_2$, every point of $\Gamma_2$ lies close to the image of $q$. Therefore $q$ is uniformly proper, by Lemma \ref{3.2} above, from $\Gamma_1$ onto a net in $\Gamma_2$.  
Hence $q$ is a quasi-isometry. This concludes the proof of the main  theorem of this subsection.

\begin{theorem}
Let $\phi$ be a uniformly proper (bijective, by definition) map from $\JJ_1 \rightarrow \JJ_2$. There exists a quasi-isometry $q$ from $\Gamma_1$ to $\Gamma_2$ which pairs the sets $\JJ_1$ and $\JJ_2$ as $\phi$ does.
\label{qipairs}
\end{theorem}

We have thus shown one aspect of {\bf relative rigidity}, viz. upgrading a uniformly proper map between $\JJ_i$'s to a quasi-isometry between $\Gamma_i$'s. In the next subsection, we shall deduce the second aspect, viz. isomorphism of $C$-complexes.

\subsection{C-Complexes}

By Theorem \ref{qipairs} we obtain a quasi-isometry $q$ from $\Gamma_1$ to $\Gamma_2$ which pairs $\JJ_1$ and $\JJ_2$ as $\phi$ does. Since $q$ is a quasi-isometry, it extends to a quasiconformal homeomorphism from $\partial G_1$ to $\partial G_2$. Also, for all $\alpha > 0$, there exists $\beta > 0$ such that $$d_1(x, J^1_i) \leq \alpha \Rightarrow d_2(q(x), \phi  (J^1_i)) \leq \beta$$ and conversely, 
$$d_2(y, J^2_i) \leq \alpha \Rightarrow d_1(x, \phi^{-1}(J^2_i)) \leq \beta$$.

In particular, $\partial q$ maps the limit set $L_i^1$ homeomorphically to the limit set of $\phi  (J^1_i))$. Hence, $\partial q$ preserves intersection patterns of limit sets.
Since $\phi$ pairs $\JJ_1$ with $\JJ_2$ as $q$ does,
summarising we get:

\begin{lemma}
The following are equivalent. \\
1) $\cap_{i=1}^k L_i^1 = \emptyset $ \\
2) $\cap_{i=1}^k \partial q(L_i^1) = \emptyset $ \\
3) $\cap_{i=1}^k \phi (L_i^1) = \emptyset $ 
\label{intpattn}
\end{lemma}

Hence by the 
definition of the $C$-complexes $C(G_1,H_1)$ and $C(G_2,H_2)$, we find that $\partial q$ induces an isomorphism of $C(G_1,H_1)$ with $C(G_2,H_2)$. We conclude:

\begin{theorem}
Let $\phi : \JJ_1 \rightarrow \JJ_2$ be a uniformly proper map. Then $\phi$ induces an isomorphism of $C(G_1,H_1)$ with $C(G_2,H_2)$.
\label{isomccx}
\end{theorem}

\noindent {\bf Note:} In Theorem \ref{qipairs} and Theorem \ref{isomccx} we start with the assumption that there exists a uniformly proper pairing of the collections $\JJ_1$ and $\JJ_2$. This can be translated to a pairing of collections of limit sets $\LL_1$ and $\LL_2$. Theorem \ref{qipairs} then says that the pairing of the $\JJ_i$'s (or $\LL_i$'s) is induced by a quasi-isometry from $\Gamma_1$ to $\Gamma_2$. Thus, the existence of a uniformly proper pairing implies the existence of a quasi-isometry between the $\Gamma_i$'s, i.e. an ambient extension (or, equivalently,  a quasiconformal homeomorphism between $\partial G_i$'s). 

Also Theorem \ref{isomccx} shows that a uniformly proper pairing induces an isomorphism of the $C$-complexes $C(G_i,H_i)$. This is reminiscent of the initial step in the proof of rigidity theorems for higher rank symmetric spaces, where Tits complexes replace $C$-complexes. 

\subsection{Cross Ratios, Annular Systems and a Dynamical Formulation}

In this subsection, we give a more intrinsic formulation of Theorems \ref{qipairs} and \ref{isomccx}. The hypothesis of these theorems is given in terms of distances between elements of $\JJ_i$. A more intrinsic way of formulating this hypothesis would be in terms of the action of $G_i$ on $\partial G_i$, $i = 1, 2$. 
In this case, the distance between $J^i_l, J^i_m$ can be approximated by the hyperbolic cross-ratio of their limit sets. This was described in detail by Bowditch \cite{bowditch-jams}. We give the relevant definitions and Theorems below and then dynamically reformulate Theorems \ref{qipairs} and \ref{isomccx}.

Let $M$ be a compactum.

\begin{defn}
An {\bf annulus} $\AAA$ is an ordered pair $(A^-,A^+)$ of disjoint closed subsets of $M$ such that $M \setminus (A^- \cup  A^+) \neq \emptyset$. An {\bf annulus system} is a collection of such annulii. If $A=(A^-,A^+)$, then $-A=(A^+,A^-)$. An annulus system is symmetric if $A\in \AAA \Rightarrow -A \in \AAA$. \\
Given a closed set $K \subset M$ and an annulus $A$, we say that $K < A$ if $K \subset int A^-$. Also, $A < K$ if $K < -A$. \\
If $A, B$ are annulii, we say that $A < B$ if $M = int A^- \cup int B^+$. \\
Fix an annulus system $\AAA$. Given closed sets $K, L \subset M$, we say that the {\bf annular cross-ratio} $(K|L)_\AAA \in \natls \cup \infty$ for the maximal number $n \in \natls$ such that we can find annulii $A_1, \cdots A_n \in \AAA$ such that $$K < A_1 < \cdots < A_n < L$$. We set $(K|L)_\AAA = \infty$ if there is no such bound.
\end{defn}

Thus $(K|L)_\AAA$ is the length of the maximal chain of nested annulii sepatrating $K, L$. For two point sets $\{ x,y \} = K$
and $\{ z,w \} = L$, we write $(K|L)_\AAA$ as $(xy|zw)_\AAA$. 

One of the crucial results of \cite{bowditch-jams} is:

\begin{theorem} (Bowditch \cite{bowditch-jams}) 
Suppose a group $G$ acts as a uniform convergence group on a perfect metrizable compactum $M$. Then there exists a symmetric $G$-invariant annulus system $\AAA$ such that if $x,y, z, w$ are distinct elements in $M$, then the theree quantities
$(xy|zw)_\AAA$, $(xz|yw)_\AAA$, $(xw|zy)_\AAA$ are all finite and at least two of them are zero. Also, if $x \neq y$, then $(x|y)_\AAA > 0$. Further, $G$ is hyperbolic, and $d_G(J(K),J(L))$ differs from $(K,L)_\AAA$ upto  bounded additive and multiplicative factors. 
\label{acr}
\end{theorem}

Combining Theorems \ref{qipairs} , \ref{isomccx} with Proposition \ref{qc=discrete} and Theorem \ref{acr}, we get the dynamical formulation we promised. Let $C_c^0(M)$ denote the collection of closed subsets of $M$ containing more than one point. (Replacing $d_{G_i}$ by cross-ratios
$(..|..)_i$ in Definition \ref{uproper} we get the corresponding notion of a map being uniformly proper with respect to the cross-ratios $(..|..)_1$ and $(..|..)_2$
in the theorem
below. Similarly for the homeomorphism $q$.)

\begin{theorem}
Let $G_1, G_2$ be uniform convergence  (hence hyperbolic) groups acting on compacta $M_1, M_2$ respectively. Also, let $\AAA_i$ (for $i = 1, 2$ ) be $G_i$-invariant annulus systems and let $(..|..)_i$ denote the corresponding annular cross-ratios.\\
 Let $H_1, H_2$ be subgroups of $G_1, G_2$ with limit sets $\Lambda_1, \Lambda_2$. Suppose that the set $\LL_i$ of translates of $\Lambda_i$ (for $i = 1, 2$) by essentially distinct elements of $H_i$ in $G_i$ forms a discrete subset  of $C_c^0(M_i)$. \\
Also assume that there exists a bijective function $\phi : \LL_1 \rightarrow \LL_2$ and  that this pairing is uniformly proper with respect to the cross-ratios $(..|..)_1$ and $(..|..)_2$. \\ Then 
\begin{enumerate}
\item $H_i$ is quasiconvex in $G_i$
\item There is a homeomorphism $q: M_1 \rightarrow M_2$ which pairs $\LL_1$ with $\LL_2$ as $\phi$ does. Further, $q$ is uniformly proper with respect to the cross-ratios $(..|..)_1$ and $(..|..)_2$ on $M_1$, $M_2$ respectively.
\item $q$ (and hence also $\phi$) induces an isomorphism of $C$-complexes $C(G_1,H_1)$ with $C(G_2,H_2)$.
\end{enumerate}
\label{dyna}
\end{theorem}

Thus from a uniformly proper map with respect to the pseudometrics on $\LL_i$'s induced by cross-ratios  we infer a quasi-isometry that is an ambient extension as also a (simplicial) isomorphism of $C$-complexes. 

\subsection{Axiomatisation, Relative Hyperbolicity} 

For classes of pairs $(X,\JJ )$, what did we really require  to ensure relative rigidity? Assume $(X,d)$ is a metric space and let the induced pseudometric on $\JJ$ be also denoted by $d$. \\
1) For all $k>0$ there exists $M \in \natls$ such that for all $x \in X$, $N_k(x)$ meets at most $M$ of the $J$'s in $\JJ$. {\it (This is a coarsening of the notion of height.)} \\
2) For all $K \in \natls$, there exists $k = k(K) > 0$ such that for all $x \in X$, $N_k(x)$ meets at least $K$ of the $J$'s in $\JJ$. {\it (This is the converse condition to (1).)} \\
3) For all $k >0, n \in \natls$ there exists $K > 0$ such that for any collection $J_1, \cdots , J_n \in \JJ$ with $d(J_i,J_j) \leq k$, there exists a ball of radius at most $K$ meeting all the $J_i$'s. \\
4) There exists $N  \in \natls$ such that for all $k > 0$ there exists $K = K(k) > 0$ such that the following holds. \\
For all $n \geq N$ and $J_1, \cdots , J_n \in \JJ$, 
 the set of points  $\{ x \in X: N_k(x) \cap J_i \neq \emptyset , i = 1\cdots n \}$ is either empty or has diameter bounded by $K$.  \\

\smallskip

Given (1)-(4), the construction of $q: X_1 \rightarrow X_2$ from a uniformly proper pairing $\phi : \JJ_1 \rightarrow \JJ_2$ goes through as in Theorem \ref{qipairs}. In short, pick $N$ from (4). From (2), pick $k = k(N)$. Now for all $x \in X_1$, consider the collection of $J$'s in $\JJ_1$ that meet $N_k(x)$. By (1) there is an upper bound $M=M(k)$ on the number of such $J$'s. Map these over to $\JJ_2$. Any two of these are at a distance of at most $m$ apart where $m$ depends on $\phi$ and $k$. From (3) choose $K = K(M,m)$ such that a ball of radius $K$ meets all these. Set $q(x)$ to be the center of such a ball. By (4), $q(x)$ is defined upto a uniformly bounded degree of discrepancy. The rest of the proof goes through as before. Hence (1)-(4) define sufficient  conditions for relative rigidity for a class of pairs $(X, \JJ )$.

\medskip

With these conditions, it is easy to extend
Theorem \ref{qipairs} to pairs $(X, \JJ )$ where $X$ is (strongly) hyperbolic relative to the collection $\JJ$. Conditions (1) and (2) are trivial. 

Condition (3) follows from ``bounded penetration" (see Farb \cite{farb-relhyp}). For any subcollection $\JJ_1 $
of $ \JJ$ with $d(J_i,J_j) \leq C_0$ (for all
$J_i, J_j \in \JJ_1$), fix
 any two 
$J_1, J_2 \in \JJ$ and a geodesic $\gamma_{12}$
of length $\leq C_0$ joining them. 
Construct an electric triangle for  triples $J_1, J_2, J_3 \in \JJ_1$ of horosphere-like sets
for arbitrary $J_3 \in \JJ_1$, such that the hyperbolic geodesics $\gamma_{13}, \gamma_{32}$ joining $J_1, J_3$
and  $J_3, J_2$ respectively have lengths
 bounded by $C_0$. 
Then $\gamma_{12}$ and $\gamma_{13}$ meet $J_1$ at a uniformly bounded distance from each other by bounded penetration. To see
this, first note that  $J_1, J_2$ can be joined by two paths, one consisting of one side of the triangle and the
other the union of the
 two remaining sides of the triangle and both paths
have electric length bounded by $2C_0$; in particular both paths are uniform quasigeodesics (with
quasigeodesic constant  depending only on $C_0$). They may be converted to quasigeodesics without backtracking 
by not increasing lengths. Thus $\gamma_{13}
\cup \gamma_{32}$ decomposes as the union of
a  quasigeodesic without backtracking 
$\gamma_{12}^{\prime}$
joining $J_1, J_2$ and (possibly) a uniformly bounded
($\leq C_0$) number of
loops of length
not longer than $2C_0$. The entire quasigeodesic 
without backtracking
$\gamma_{12}^{\prime}$ lies near $\gamma_{12}$ for all 
$J_3 \in \JJ_1$. The same holds for the loops of bounded
length (since they in turn may be regarded as
uniform quasigeodesics 
without backtracking starting and ending at the same point.) In particular $J_3$ lies at a uniformly bounded
distance $D_0$ from $\gamma_{12}$.
Since $\gamma_{12}$ has length bounded by $C_0$,
and $J_3$ may be chosen arbitrarily satisfying the hypothesis of (3) above, it follows that for any $x\in \gamma_{12}$, $d(x, J_3) \leq (C_0 + D_0)$
for all $J_3 \in \JJ_1$.  Condition (3) follows. (Results closely
related to the proof of Condition (3) here occur
 as Lemma 3.11 of \cite{brahma-ibdd}
and Prop. 8.6 of \cite{hruska-wise}.)

Condition (4) follows from the fact that for a pair of distinct $J_i, J_j$, 
$N_k(J_i) \cap N_k(J_j)$ is either empty or has diameter bounded by some $C(k)$.

We have thus shown:

\begin{theorem}
Let $X_i$ be (strongly) hyperbolic relative to collections $\JJ_i$ ($i = 1, 2$).
Let $\phi$ be a uniformly proper (bijective, by definition) map from $\JJ_1 \rightarrow \JJ_2$. There exists a quasi-isometry $q$ from $X_1$ to $X_2$ which pairs the sets $\JJ_1$ and $\JJ_2$ as $\phi$ does.
\label{qipairs-relhyp}
\end{theorem}

By work of Hruska and Kleiner \cite{hruska-kleiner}, CAT(0) spaces with isolated flats are (strongly) hyperbolic relative to maximal flats. Hence we have from Theorem \ref{qipairs-relhyp} above:

\begin{cor}
Let $X_i$ be CAT(0) spaces with isolated flats and let $\JJ_i$ denote the collections of maximal flats ($i = 1, 2$).
Let $\phi$ be a uniformly proper (bijective, by definition) map from $\JJ_1 \rightarrow \JJ_2$. There exists a quasi-isometry $q$ from $X_1$ to $X_2$ which pairs the sets $\JJ_1$ and $\JJ_2$ as $\phi$ does.
\label{qipairs-cat0}
\end{cor}

\subsection{Symmetric Spaces of Higher Rank}
We now consider CAT(0) spaces which are at the other end of the spectrum. Let $M$ be a compact locally symmetric space and $T$ a totally geodesic torus with rank = rank$(M)$. Take $X = \til{M}$ and $\JJ$ to be the lifts of $T$ to $\til{M}$. As these are all equivariant examples
(i.e. $\mathcal{J}$ is invariant under a cocompact group action), it is enough to check (1)-(4) at a point.

(1) and (2) are clear. To prove condition (4), we consider 
$\cap_i N_k(F_i)$ and it is easy to bound from below the $N$ appearing in Condition (4) (Section 3.4) in terms of the size of the Weyl group and rank. In that case, 
$\cap_i N_k(F_i)$ has bounded diameter or is empty. 

 Finally, to prove (3), we proceed as in Lemma \ref{3.1}. As in Lemma \ref{3.1} we assume by induction that any $k$ flats $\{ F_1 , \cdots , F_k \}$ that "coarsely pairwise intersect at scale $D$" 
(i.e. $N_D(F_i) \cap N_D(F_j) \neq \emptyset$ ) intersect coarsely (i.e. $\cap_{i= 1 \cdots k} N_{D^{\prime}} (F_i)  \neq \emptyset$  for some 
$D^{\prime} = D^{\prime} (D,k)$). To get to the inductive step, we suppose that for $i = k+1$, we have collections of worse and worse counterexamples. 
Consider a maximal collection $\FF = \{ F_1 , \cdots , F_k \}$
of maximal flats whose "coarse intersection at scale $D$" 
$\cap_i N_D(F_i) = F$
is non-null. Translate the collection by a group element so that a fixed point $0$ (thought of as the origin)
lies on the intersection $F$. 
 Now take a sequence of maximal flats $F^j$ whose $D$-neighborhoods
$N_D(F^j)$ intersect each $N_D(F_i)$, but  $d_j = d(F^j,F) = d(0,F) \geq j$. 
We scale the metric on $(X,d)$ by a factor of $d_j$ to obtain a sequence of metric spaces $(X, \frac{d}{d_j} )$ converging (via a non-principal ultrafilter)
to a Euclidean building
$X^\infty$ (this fact is due to Kleiner and Leeb \cite{kleiner-leeb}, but 
we shall only mildly need the exact nature of $X_\infty$). $F_i$'s converge to flats $F_i^\infty \subset X^\infty$ and $F^j$'s converge to a flat
$G^\infty \subset X^\infty$. Then the collection $\GG = F_i^\infty , 
G^\infty$ satisfy the following: \\
$(P1)$ Each element of $\GG$ is a  flat in $X^\infty$ \\
$(P2)$ By induction, the intersection of any $i$ elements of $\GG$ is non-empty and convex for $i \leq k$ \\
$(P3)$ The intersection of all the $(k+1) $ elements of $\GG$ is empty. \\

\smallskip

Consider the subcomplex $K = G^\infty \bigcup_i F_i^\infty  
$ of $X^\infty$. Then $K$ is a union of $r$-flats, where $r = {\rm rank}
(X)$. In particular, the homology groups $H_n (K) = 0$ for $n > r$.
On the other hand, if we consider the nerve of the covering of $K$ by the sets $G^\infty , F_i^\infty  
$, then using the  properties $(P1), (P2), (P3)$
to compute Cech homology groups, we conclude that
$K$ has the same homology groups as the boundary of a $k$-simplex.
In particular, $H_k(K) = \Bbb{Z}$. For $k > r$ this is a contradiction,
finally proving Condition (3). Thus we conclude: \\

\begin{theorem}
Let $X_i$ be symmetric spaces of non-positive curvature, and  $\JJ_i$ be equivariant collections of lifts of a maximal torus in a compact locally symmetric space modeled on $X_i$ ($i = 1, 2$).
Let $\phi$ be a uniformly proper (bijective, by definition) map from $\JJ_1 \rightarrow \JJ_2$. There exists a quasi-isometry $q$ from $X_1$ to $X_2$ which pairs the sets $\JJ_1$ and $\JJ_2$ as $\phi$ does.
\label{qipairs-ss}
\end{theorem}

Combining Theorem \ref{qipairs-ss} with the quasi-isometric rigidity
theorem of Kleiner-Leeb \cite{kleiner-leeb} and Eskin-Farb \cite{eskin-farb-jams}
we can upgrade the quasi-isometry $q$ to an isometry $i$.

\begin{cor}
Let $X_i$ be symmetric spaces of non-positive curvature, and  $\JJ_i$ be equivariant collections of lifts of a maximal torus in a compact locally symmetric space modeled on $X_i$ ($i = 1, 2$).
Let $\phi$ be a uniformly proper (bijective, by definition) map from $\JJ_1 \rightarrow \JJ_2$. There exists an isometry $i$ from $X_1$ to $X_2$ which pairs the sets $\JJ_1$ and $\JJ_2$ as $\phi$ does.
\label{ipairs-ss}
\end{cor}

\begin{rmk}
The technique of using asymptotic cones and the nerve of the covering
by flats can be generalised easily to equivariant flats of arbitrary
(not necessarily maximal) rank.
\end{rmk}

We conclude this paper with two (related) questions: \\
{\bf Question 1:} In analogy with a Theorem of Ivanov, Korkmaz, Luo (see for instance \cite{luo-topology} ), regarding the automorphism group of the curve complex, we ask: \\ If the C-Complex $C(G,H)$ of a pair $(G,H)$ (for $G$ a hyperbolic group and $H$ a quasiconvex subgroup) is connected, is the automorphism group of  
$C(G,H)$ commensurable with $G$? 

\medskip

\noindent {\bf Question 2:} Consider the
 pair $(G,H)$, with
$G$ a hyperbolic group and $H$ a quasiconvex subgroup.
Let $(\JJ , d)$ be the collection  of joins as in Lemma \ref{3.1} with the induced pseudometric. For a uniformly proper map $\phi$ from $(\JJ , d)$ to itself, is there an isometry pairing the elements of $\JJ$ as $\phi$? We have proved in Theorem \ref{qipairs} that a quasi-isometry $q$ exists pairing the $\JJ$ as $\phi$ does. The question is whether $q$ may be upgraded to an isometry, or better, to an element of $G$? This question is related to the notion of {\it pattern rigidity} introduced by Mosher, Sageev and Whyte in \cite{msw2}.

\bibliography{qcrh}

\begin{thebibliography}{GMRS97}

\bibitem[BCG98]{BCG}
G.~Besson, G.~Courtois, and S.~Gallot.
\newblock {A real Schwarz lemma and some applications}.
\newblock {\em Rend. Mat. Appl., VII. Ser. 18, No.2}, pages 381--410, 1998.

\bibitem[BDM05]{behrstock-drutu-mosher}
J.~Behrstock, C.~Drutu, and L.~Mosher.
\newblock Thick metric spaces, relative hyperbolicity, and quasi-isometric
  rigidity.
\newblock {\em preprint, arXiv:math.GT/0512592}, 2005.

\bibitem[Bow97]{bowditch-relhyp}
B.~H. Bowditch.
\newblock Relatively hyperbolic groups.
\newblock {\em preprint, Southampton}, 1997.

\bibitem[Bow98]{bowditch-jams}
B.~H. Bowditch.
\newblock A topological characterization of hyperbolic groups.
\newblock {\em J. A. M. S. 11}, pages 643--667, 1998.

\bibitem[DS05]{drutu-sapir}
C.~Drutu and M.~Sapir.
\newblock Tree-graded spaces and asymptotic cones of groups.
\newblock {\em Topology 44}, pages 959--1058, 2005.

\bibitem[EF 3]{eskin-farb-jams}
A.~Eskin and B.~Farb.
\newblock {Quasi-flats and rigidity in higher rank symmetric spaces. }.
\newblock {\em J. Amer. Math. Soc., 10}, pages pp. 653--692., (1997),no. 3.

\bibitem[Far98]{farb-relhyp}
B.~Farb.
\newblock Relatively hyperbolic groups.
\newblock {\em Geom. Funct. Anal. 8}, pages 810--840, 1998.

\bibitem[FM98]{farb-mosher-bs1}
B.~Farb and L.~Mosher.
\newblock {A rigidity theorem for the solvable Baumslag-Solitar groups}.
\newblock {\em Inventiones Math. , Vol. 131, No.2}, pages 419--451, 1998.

\bibitem[FM00]{farb-mosher-acta}
B.~Farb and L.~Mosher.
\newblock On the asymptotic geometry of abelian-by-cyclic groups.
\newblock {\em Acta Math. , Vol. 184, No.2}, pages 145--202, (2000).

\bibitem[FM 3]{farb-mosher-bs2}
B.~Farb and L.~Mosher.
\newblock {Quasi-isometric rigidity for the solvable Baumslag-Solitar groups,
  II}.
\newblock {\em Inventiones Math. , Vol. 137}, pages 273--296, (1999), No. 3.

\bibitem[FS96]{farb-schwarz}
B.~Farb and R.E. Schwartz.
\newblock {The large-scale geometry of Hilbert modular groups}.
\newblock {\em J. Differential Geom. 44(3)}, page 435–478, 1996.

\bibitem[GdlH90]{GhH}
E.~Ghys and P.~de~la Harpe(eds.).
\newblock Sur les groupes hyperboliques d'apres {M}ikhael {G}romov.
\newblock {\em Progress in Math. vol 83, Birkhauser, Boston Ma.}, 1990.

\bibitem[GMRS97]{GMRS}
R.~Gitik, M.~Mitra, E.~Rips, and M.~Sageev.
\newblock Widths of {S}ubgroups.
\newblock {\em Trans. AMS}, pages 321--329, Jan. '97.

\bibitem[Gro85]{gromov-hypgps}
M.~Gromov.
\newblock Hyperbolic {G}roups.
\newblock {\em in Essays in Group Theory, ed. Gersten, MSRI Publ.,vol.8,
  Springer Verlag}, pages 75--263, 1985.

\bibitem[HK04]{hruska-kleiner}
G.~Christopher Hruska and Bruce Kleiner.
\newblock {Hadamard spaces with isolated flats}.
\newblock {\em {arxiv:math.GR/0411232}}, preprint 2004.

\bibitem[HW06]{hruska-wise}
G.~Christopher Hruska and D.~T. Wise.
\newblock {Packing Subgroups in Relatively Hyperbolic Groups}.
\newblock {\em {arxiv:math.GR/0609369}}, preprint 2006.

\bibitem[KL97a]{kap-leeb}
M.~Kapovich and B.~Leeb.
\newblock {Quasi-isometries preserve the geometric decomposition of Haken
  manifolds}.
\newblock {\em Inventiones Math., Vol. 128, F. 2}, pages 393--416, (1997).

\bibitem[KL97b]{kleiner-leeb}
B.~Kleiner and B.~Leeb.
\newblock Rigidity of quasi-isometries for symmetric spaces and euclidean
  buildings.
\newblock {\em IHES Publ. Math. 86}, page 115–197, (1997).

\bibitem[Luo00]{luo-topology}
F.~Luo.
\newblock Automorphisms of the complex of curves.
\newblock {\em Topology 39}, pages 283--298., (2000).

\bibitem[Mit98]{mitra-trees}
Mahan Mitra.
\newblock Cannon-{T}hurston {M}aps for {T}rees of {H}yperbolic {M}etric
  {S}paces.
\newblock {\em Jour. Diff. Geom.48}, pages 135--164, 1998.

\bibitem[Mit04]{mitra-ht}
Mahan Mitra.
\newblock Height in {S}plittings of {H}yperbolic {G}roups.
\newblock {\em Proc. Indian Acad. of Sciences, v. 114, no.1}, pages 39--54,
  Feb. 2004.

\bibitem[Mj05]{brahma-ibdd}
Mahan Mj.
\newblock {Cannon-Thurston Maps, i-bounded Geometry and a Theorem of McMullen}.
\newblock {\em preprint, arXiv:math.GT/0511041}, 2005.

\bibitem[Mos73]{mostow-book}
G.~D. Mostow.
\newblock {Strong Rigidity of Locally Symmetric Spaces}.
\newblock {\em Princeton University Press}, 1973.

\bibitem[MSW03]{msw1}
L.~Mosher, M.~Sageev, and K.~Whyte.
\newblock {Quasi-actions on trees I.Bounded valence}.
\newblock {\em Annals of Mathematics, 158, arXiv:math.GR/0010136}, page
  115–164, (2003).

\bibitem[MSW04]{msw2}
L.~Mosher, M.~Sageev, and K.~Whyte.
\newblock {Quasi-actions on trees II.Finite depth Bass-Serre trees}.
\newblock {\em arXiv:math.GR/0405237}, (2004).

\bibitem[NR03]{niblo-reeves}
G.~Niblo and L.~Reeves.
\newblock Coxeter groups act on cat(0) cube complexes.
\newblock {\em J. Group Theory 6}, pages 399--413, (2003).

\bibitem[Pau96]{paulin-bdy}
F.~Paulin.
\newblock Un groupe hyperbolique est déterminée par son bord.
\newblock {\em J. London Math. Soc. 54}, pages 50--74, 1996.

\bibitem[Roe95]{roe-cbms}
J.~Roe.
\newblock {Index Theory, Coarse Geometry, and Topology of Manifolds}.
\newblock {\em CBMS series, AMS}, 1995.

\bibitem[Sag95]{sageev-th}
M.~Sageev.
\newblock Ends of group pairs and non-positively curved cube complexes.
\newblock {\em Proc. London Math. Soc. 71}, pages 585--617, 1995.

\bibitem[Sch97]{schwarz-inv}
R.~E. Schwarz.
\newblock Symmetric patterns of geodesics and automorphisms of surface groups.
\newblock {\em Invent. math. Vol. 128, No. 1}, page 177–199, 1997.

\bibitem[Sho91]{short}
H.~Short.
\newblock Quasiconvexity and a theorem of {H}owson's.
\newblock {\em Group {T}heory from a {G}eometrical {V}iewpoint ({E}. {G}hys,
  {A}. {H}aefliger, {A}. {V}erjovsky eds.)}, 1991.

\end{thebibliography}
\bibliographystyle{alpha}

\end{document}